\documentclass[review,1p]{elsarticle}
\usepackage{graphicx,amssymb,mathrsfs,amsmath,amsfonts,dsfont,multirow}
\usepackage{lineno}
\usepackage{subfigure}

\modulolinenumbers[12]

\newtheorem{thm}{Theorem}[section]

\newdefinition{rmk}{Remark}
\newtheorem{pro}{Proposition}[section]
\newproof{pf}{Proof}
\newproof{pot}{Proof of Theorem \ref{thm2}}

\numberwithin{equation}{section}

\newtheorem{defn}{Definition}[section]
\newtheorem{prop}{Proposition}[section]

\journal{Journal of   }

\begin{document}

\begin{frontmatter}

\title{A Nonconvex Nonsmooth Regularization Method for Compressed Sensing and Low-Rank Matrix Completion}



\author[]{Zhuo-Xu Cui}
\ead{zhuoxucui@whu.edu.cn}
\author[]{Qibin Fan\corref{mycorrespondingauthor}}

\cortext[mycorrespondingauthor]{Corresponding author}
\ead{qbfan@whu.edu.cn}

\address{School of Mathematics and Statistics, Wuhan University, Wuhan 430072, China}

\begin{abstract}
In this paper, nonconvex and nonsmooth models for compressed sensing (CS) and low-rank matrix completion (MC) is studied. The problem is formulated as a nonconvex regularized leat square optimization problems, in which the $\ell_{0}$-norm and the rank function are replaced by $\ell_{1}$-norm and nuclear norm, and adding a nonconvex penalty function respectively. An alternating minimization scheme is developed, and the existence of a subsequence, which generate by the alternating algorithm that converges to a critical point, is proved. The NSP, RIP, and RIP condition for stable recovery guarantees also be analysed for the nonconvex regularized CS and MC problems respectively. Finally, the performance of the proposed method is demonstrated through experimental results.
\end{abstract}


\end{frontmatter}

\linenumbers

\section{ Introduction}
The compressed sensing (CS) problem is to recover an unknown vector from a small amount of observations. It's possible to exactly reconstruct it with high probability if the vector is sparse. The mathematical formula reads:
$$\min_{x}\{\|x\|_{0}~:~Ax=y\},\eqno (1.1)$$
where  $x\in\mathbb{R}^{p}$, $y\in\mathbb{R}^{n}$ with $n\ll p$, $A\in\mathbb{R}^{n\times p}$ is a measurement ensembles \cite{T. Tao{2005},T. Tao{2006}a,T. Tao{2006}b,T. Tao{2006}c,M.B. Wakin{2008}}. The matrix completion (MC) problem is to recover a low-rank matrix from a small amount of observations:
$$\min_{X}\{rank(X)~:~X_{ij}=Y_{ij},~(i,j)\in\Omega \},\eqno (1.2)$$
where $X\in\mathbb{R}^{n_{1}\times n_{2}}$, $\Omega$ is a given set of index pairs $(i,j)$ \cite{B. Recht{2009},T. Tao{2010}}.

Due to the minimization of $\ell_{0}$-norm and $rank$ function, these problems (1.1), (1.2) are NP-hard problem in general, In some sense, $\ell_{1}$-norm and nuclear norm are the tightest convex relaxation of these nonconvex functions, respectively. The nuclear norm of $X$ define as $\|X\|_{*}=\sum_{i=1}^{m}\sigma_{i}(X)$, where $\sigma_{i}$ is the $i$ largest singular value of $X$ and $m$ is the number of singular value. Therefore, the problem (1.1) and (1.2) can be relaxed into:
 $$\min_{x\in\mathbb{R}^{p}}\{\|x\|_{1}~:~Ax=y\},\eqno (1.3)$$
$$\min_{X\in\mathbb{R}^{n_{1}\times n_{2}}}\{\|X\|_{*}~:~X_{ij}=Y_{ij},~(i,j)\in\Omega \},\eqno (1.4)$$
and the problem (1.3) and (1.4) is equivalent to (1.1) and (1.2) respectively under certain incoherence conditions \cite{S. Foucart{2013}}.
However, the solution of (1.3) and (1.4) is usually suboptimal to the original problem (1.1) and (1.2), the $\ell_{1}$-norm minimization problem may yield the vector with lower sparse rate than the real one, and can't recover a sparse target with minimum measurements. Another limitation of the $\ell_{1}$-norm minimization is its bias caused by shrinking all the element toward zero simultaneously \cite{Y. Hu{2013}}, the nuclear norm of a matrix is the $\ell_{1}$-norm of it's singular value vector, so it also have these limitations.

Since the $\ell_{1}$-norm may not be approximated $\ell_{0}$-norm well, in CS recovery problems, many known nonconvex surrogates of $\ell_{0}$-norm have been proposed, include $\ell_{p}$-norm($0<p<1$) \cite{L.L.E. Frank{1993}}, Smoothly Clipped Absolute Deviation (SCAD) \cite{J. Fan{2001}}, Minimax Concave Penalty (MCP) \cite{C.H. Zhang{2010}}, Exponential Type Penalty (ETP) \cite{C. Gao{2011}}, etc. Recently, some of these method have been extended to low-rank matrix restoration and have well performance.

Because of the limitation of (1.3) and (1.4), we augment them by adding a nonconvex and nonsmooth term $\beta\Phi(x)$ and $\beta\hat{\Phi}(X)$, respectively, where $\beta$ is a positive scalar,  $$\Phi(x)=\sum^{p}_{i=1}\varphi(x_{i}),~\hat{\Phi}(X)=\sum_{i=1}^{m}\varphi(\sigma_{i}(X)),\eqno (1.5)$$ where $\varphi(t)=\frac{\alpha|t|}{1+\alpha|t|}$ \cite{Nikolova M{2012},X. Chen{2010}}, $\sigma_{i}$ is the $i$ largest singular value of $X$ and $m$ is the number of singular value. The augmented model for (1.3) and (1.4) are
$$\min_{x\in\mathbb{R}^{p}}\{\|x\|_{1}+\beta\Phi(x):Ax=y\},\eqno (1.6)$$
$$\min_{X\in\mathbb{R}^{n_{1}\times n_{2}}}\{\|X\|_{*}+\beta\hat{\Phi}(X):\mathcal{A}X=b\},\eqno (1.7)$$
which can be solved by introducing a auxiliary variable and using alternating minimization scheme \cite{Jin X{2015}}. In (1.7), $\mathcal{A}$ is a linear operator, if we choose $\mathcal{A}$ as a componentwise projection, it become the matrix completion problem.
The solution to (1.6) and (1.7) is also a solution to (1.3) and (1.4) as long as $\beta$ is sufficiently small, and $\beta$ controls the tradeoff between $\ell_{1}$-norm term and nonconvex term. For recovering a sparse vector and a low-rank matrix, the choose of the suitable $\beta$ should obey follow formula
$$\beta\leq\frac{1}{20\alpha}.\eqno (1.8)$$
In general, we choose $\alpha=0.5$, so $\beta\leq0.1$.

One can observe that $\Phi(x)$ convergence to $\|x\|_{0}$ and $C\|x\|_{1}$, as $\alpha\rightarrow\infty$ and $\alpha\rightarrow 0$ respectively, where $C$ is a large scaler. It has been show in \cite{Mila Nikolova{2010}} that $\varphi$ satisfies: (1) $\varphi$ is continuous (Lipschitz function), symmetric on $(-\infty,\infty)$, $C^{2}$ on $(0,\infty)$ and $\varphi(0)=0$ is a strict minimum; (2) $\varphi'(0^{+})>0$ and $\varphi'(t)\geq0$ for all $t>0$; (3) $\varphi''$ is increasing on $(0,\infty)$ with $\varphi''(t)<0$ and $\lim_{t\rightarrow\infty}\varphi''(t)=0$, which implies that our augment regularizers to be a good promoted penalty function, and the augment term have some properties as follows:\\
(1) $\alpha>0$, $\hat{\Phi}(X)\geq0(\Phi(x)\geq0)$, with equality hold if only if $X=0(x=0)$;\\
(2) $\hat{\Phi}(X)(\Phi(x))$ is a decreasing function of $\alpha$, and $\lim_{\alpha\rightarrow\infty}\hat{\Phi}(X)=rank(X)(\lim_{\alpha\rightarrow\infty}\Phi(x)=\|x\|_{0})$;\\
(3)  $\hat{\Phi}(X)$ is unitarily invariant, that is $\hat{\Phi}(UXV^{*})=\hat{\Phi}(X)$ whenever $U\in\mathbb{R}^{n_{1}\times n_{1}}$ and $V\in\mathbb{R}^{n_{2}\times n_{2}}$ are orthogonal matrix.

This paper also shows the recovery guarantees for augment model of compressed sensing and low-rank matrix completion respectively, the results are given based on varieties of properties of matrix $A$ and linear operator $\mathcal{A}$ including the null-space property (NSP), the restricted isometry property (RIP), at last, the RIP condition for stable recovery are given.

The rest of this paper is organized as follows. In Sect. 2, we firstly give the augmented model, and introduce the nonconvex and nonsmooth penalty function for low-rank matrix completion and sparse vector recovery. Then, we use the alternating minimization scheme for solving the proposed problem and give the convergence result of the proposed method. In Sect. 3, we shows the recovery guarantees for augmented model of compressed sensing and low-rank matrix completion respectively, include NSP, RIP, and so on.
In Sect. 4, some numerical experiment results of our augment model have been showed on simulated and real data. Finally, some conclusions are summarized in Sect. 5.

\section{Algorithm and Convergence Analysis}
In this section, we propose an alternating minimization scheme for solving (1.6) and (1.7). We begin with introducing an auxiliary variable, and obtain a new cost function, then we decompose the cost function into two subproblems, soft-thresholding operator has been used to solve subproblem one and Quasi-Newton¡¯s method has been used to solve subproblem two. Finally, we give the algorithm for solving (2.5) and show its convergence.

Firstly, we consider the variant of (1.6) and (1.7) are
$$\min_{x\in\mathbb{R}^{p}}\{\|x\|_{1}+\beta\Phi(x):\|Ax-y\|_{2}\leq\epsilon\},\eqno (2.1)$$
$$\min_{X\in\mathbb{R}^{n_{1}\times n_{2}}}\{\|X\|_{*}+\beta\hat{\Phi}(X):\|\mathcal{A}X-b\|_{2}\leq\epsilon\},\eqno (2.2)$$
where $\epsilon\geq0$ admits the possible noise in the measurement. The equivalent Lagrangian form:
$$\min_{x\in\mathbb{R}^{p}}:~\frac{1}{2}\|Ax-y\|_{2}^{2}+\lambda(\|x\|_{1}+\beta\Phi(x)),\eqno (2.3)$$
$$\min_{X\in\mathbb{R}^{n_{1}\times n_{2}}}:~\frac{1}{2}\|\mathcal{A}X-b\|_{2}^{2}+\lambda(\|X\|_{*}+\beta\hat{\Phi}(X)),\eqno (2.4)$$
where $\lambda$ is the regularization parameter which controls the tradeoff between data fitting term and the regularization term. Next, we mainly introduce the low-rank matrix completion problems, and it is fairly easy to extended the result to sparse vector recovery.

Firstly, by introducing an auxiliary variable $W\in\mathbb{R}^{n_{1}\times n_{2}}$, cost function (2.4) can be approximately transformed into
$$\varepsilon(X,W)=\frac{1}{2}\|\mathcal{A}X-b\|_{2}^{2}+\lambda(1+\alpha\cdot\beta)\|W\|_{*}+\lambda\cdot\beta\hat{\Psi}(X)+\frac{\rho}{2}\|X-W\|^{2}_{F},\eqno (2.5)$$
where $\hat{\Psi}(X)=\hat{\Phi}(X)-\alpha\cdot\|X\|_{*}$, and there exists Gateaux derivatives of $\hat{\Psi}(X)$ at $X$, however, the Gateaux derivatives of $\hat{\Phi}(X)$ is not always exist.

Given $(W^{(s-1)},X^{(s-1)})$, the iteration scheme of problem (2.5) can be described as follows:
$$~~~~W^{(s)}\in\arg\min_{W\in\mathbb{R}^{n_{1}\times n_{2}}}\varepsilon(W,X^{(s-1)});\eqno (2.6)$$
$$X^{(s)}\in\arg\min_{X\in\mathbb{R}^{n_{1}\times n_{2}}}\varepsilon(W^{(s)},X),\eqno (2.7)$$
where $\arg\min$ denotes the minimal set to an optimization problem.
It's easy to know that the W-subproblem (2.6) can formulated as
$$~~W^{(s)}\in \arg\min_{W\in\mathbb{R}^{n_{1}\times n_{2}}}\frac{\rho}{2}\|X-W\|^{2}_{F}+\lambda(1+\alpha\cdot\beta)\|W\|_{*}$$
$$= \arg\min_{W\in\mathbb{R}^{n_{1}\times n_{2}}}\frac{1}{2}\|X-W\|^{2}_{F}+\tau\|W\|_{*},\eqno (2.8)$$
where $\tau=\frac{\lambda(1+\alpha\cdot\beta)}{\rho}$, according to \cite{J.F.Cai{2010}}, it's easy to show the solution of (2.8) as
$$W^{(s)}=\mathcal{D}_{\tau}(X^{(s-1)}),\eqno (2.9)$$
where $\mathcal{D}_{\tau}$ is the soft-thresholding operator, $\mathcal{D}_{\tau}=U\mathcal{D}_{\tau}(\Sigma)V^{*}$, $\mathcal{D}_{\tau}(\Sigma)=diag(\{\sigma_{i}-\tau\}_{+})$, $t_{+}$ is the positive part of $t$, namely, $t_{+}=\max(0,t)$ and $X=U\Sigma V^{*}$ is the singular value decomposition (SVD) of matrix $X$.

The X-subproblem (2.7) can be formulated as follows
$$X^{(s)}=\arg\min_{X\in\mathbb{R}^{n_{1}\times n_{2}}}\frac{1}{2}\|\mathcal{A}X-b\|^{2}_{2}+\lambda\cdot\beta\hat{\Psi}(X)+\frac{\rho}{2}\|X-W^{(s)}\|^{2}_{F},\eqno (2.10)$$
we could use Quasi-Newton's method to solve this optimization problem
$$(\mathcal{A}^{*}\mathcal{A}+\rho \mathcal{I})\Delta X=\mathcal{A}^{*}(b-\mathcal{A}X^{(s-1)})-\lambda\cdot\beta D_{X}\hat{\Psi}(X^{(s-1)})+\rho(W^{(s)}-X^{(s-1)}),\eqno (2.11)$$
$$X^{(s)}=X^{(s-1)}+\Delta X,\eqno (2.12)$$
where $\mathcal{I}$ is an identity operator, and $\mathcal{A}^{*}$ is the adjoint of $\mathcal{A}$. In order to get $\Delta X$, we could use conjugate gradient method for solving this linear system (2.11).
\\

\begin{pro}
The Gateaux derivatives of $\hat{\Psi}(X)$ is
$$D_{X}\hat{\Psi}(X)=U^{*}\Lambda V,\eqno (2.13)$$
where $\Lambda=Diag(\frac{\partial\psi(\sigma_{1})}{\partial\sigma_{1}},\ldots,\frac{\partial\psi(\sigma_{m})}{\partial\sigma_{m}})_{n_{1}\times n_{2}}$, $\psi(t)=\varphi(t)-\alpha|t|$ and $U$, $V$ are  unitary matrices which consist of left-singular vectors and right-singular vectors.
\end{pro}
\begin{pf}
$\varphi$ is a nonconvex and nonsmooth function, and $\varphi(t)=\alpha|t|+\psi(t)$, $\psi\in C^{2}$.
$D_{\sigma}\hat{\Psi}(\sigma)=Diag(\frac{\partial\psi(\sigma_{1})}{\partial\sigma_{1}},\ldots,\frac{\partial\psi(\sigma_{s})}{\partial\sigma_{s}})_{n_{1}\times n_{2}}$,
$\Sigma(X)=U^{*}XV$, $U$ and $V$ are unitary matrices which consist of left-singular vectors and right-singular vectors, and $\Sigma(X)\in\mathbb{R}^{n_{2}\times n_{2}}\rightarrow\mathbb{R}^{n_{2}\times n_{2}}$, we have $D_{X}\Sigma(X)\in\mathbb{R}^{n_{2}\times n_{2}}\rightarrow L(\mathbb{R}^{n_{2}\times n_{2}},\mathbb{R}^{n_{2}\times n_{2}})$, $\langle D_{X}\Sigma(X),H\rangle=U^{*}HV$, where $H\in\mathbb{R}^{n_{2}\times n_{2}}$ is an arbitrary matrix. By chain rule of Gateaux derivatives, we have $D_{X}\hat{\Psi}(X)=U^{*}Diag(\frac{\partial\psi(\sigma_{1})}{\partial\sigma_{1}},\ldots,\frac{\partial\psi(\sigma_{s})}{\partial\sigma_{s}})_{n_{1}\times n_{2}} V$.
\end{pf}

Based on the analysis above, we give a basic framework of the alternating minimization scheme for solving our nonconvex augmented model of low-rank completion problem as follows:

\begin{tabular}{l}
 \hline
Algorithm to Solve The Minimum Value of (2.5) \\
\hline
Step 1: Initialize $X^{(0)}$ and $s=1$;\\
Step 2: Update $X$ and $W$ until the convergence\\
~~~~~~~~W-step:\\
~~~~~~~~~~~$W^{(s)}=\arg\min_{W\in\mathbb{R}^{n_{1}\times n_{2}}}\varepsilon(W,X^{(s-1)})$,\\
~~~~~~~~~~~$W^{(s)}=\mathcal{D}_{\tau}(X^{(s-1)})$ and $\tau=\frac{\lambda(1+\alpha)}{\rho}$.\\
~~~~~~~~X-step:\\
~~~~~~~~~~~$X^{(s)}=X^{(s-1)}+\Delta X$, where\\
~~~~~~~~~~~$(\mathcal{A}^{*}\mathcal{A}+\rho \mathcal{I})\Delta X=-D_{X} \varepsilon(X,W^{(s)})$,\\
~~~~~~~where, $D_{X}\varepsilon(X,W^{(s)})$ is the Gateaux derivatives at X.\\
(Here the iteration index is the superscript $s$.)\\
\hline
\end{tabular}\\

\begin{prop}

(1) For all $s\geq1$, there exist a $\gamma$ such that
$$\varepsilon(W^{(s)},X^{(s)})+\gamma(\|W^{(s)}-W^{(s-1)}\|_{F}^{2}+\|X^{(s)}-X^{(s-1)}\|_{F}^{2})\leq\varepsilon(W^{(s-1)},X^{(s-1)}),\eqno (2.14)$$
hence, $\varepsilon(W^{(s)},X^{(s)})$ dose not increase.\\
(2) $$\sum_{s=1}^{\infty}(\|W^{(s)}-W^{(s-1)}\|_{F}^{2}+\|X^{(s)}-X^{(s-1)}\|_{F}^{2})<+\infty.\eqno (2.15)$$ \cite{Attouch H{2008},Attouch H{2005}}.
\end{prop}

\begin{thm}
Let $\{(W^{(s)},X^{(s)})\}$ be a sequence generated by our algorithm, then there exists a subsequence of $\{(W^{(s)},X^{(s)})\}$ such that it converges to a critical point.
\end{thm}
\begin{pf}
According to (2.8), we first obtain
$$0\in\partial_{W}\varepsilon(W^{(s)},X^{(s-1)})=\lambda(1+\alpha\cdot\beta)\partial\|W\|_{*}+\rho(W^{(s)}-X^{(s-1)}),\eqno (2.16)$$
and we have
$$~~~~~\partial_{W}\varepsilon(W^{(s)},X^{(s)})=\lambda(1+\alpha\cdot\beta)\partial\|W\|_{*}+\rho(W^{(s)}-X^{(s)}).\eqno (2.17)$$
According to (2.12), we obtain
$$-(\mathcal{A}^{*}\mathcal{A}+\rho \mathcal{I})(X^{(s)}-X^{(s-1)})=D_{X} \varepsilon(X^{(s-1)},W^{(s)})$$
$$=\mathcal{A}^{*}(Y-\mathcal{A}X^{(s-1)})-\lambda\cdot\beta D_{X}\hat{\Psi}(X^{(s-1)})+\rho(W^{(s)}-X^{(s-1)}),\eqno (2.18)$$
and we have
$$D_{X} \varepsilon(X^{(s)},W^{(s)})=\mathcal{A}^{*}(\mathcal{A}X^{(s)}-Y)+\lambda\cdot\beta D_{X}\hat{\Psi}(X^{(s)})+\rho(X^{(s)}-W^{(s)}).\eqno (2.19)$$
With (2.16), (2.17) and (2.18), (2.19), we obtain
$$\partial_{W}\varepsilon(W^{(s)},X^{(s)})=\rho(X^{(s-1)}-X^{(s)}),\eqno (2.20)$$
and
$$D_{X} \varepsilon(X^{(s)},W^{(s)})=\lambda\cdot\beta (D_{X}\hat{\Psi}(X^{(s)})-D_{X}\hat{\Psi}(X^{(s-1)})).\eqno (2.21)$$
Suppose there exist a bounded subsequence $\{(W^{(s')},X^{(s')})\}$, by using (2.15) we have
$$\lim_{s\rightarrow+\infty}\{(W^{(s')},X^{(s')})\}-\{(W^{(s'-1)},X^{(s'-1)})\}=0,\eqno (2.22)$$
and $D_{X}\hat{\Psi}(X)$ is a continuous function on bounded subsets, then,
$$\{(W^{*},X^{*})\}=\lim_{s\rightarrow+\infty}\{(W^{(s')},X^{(s')})\},\eqno (2.23)$$ is a critical point.
\end{pf}

\section{Recovery Guarantees}
In this section, we established recovery guarantees for our augmented models (1.6) and extends these result to matrix recovery models (1.7). The result for (1.6) and (1.7) are given based on varieties of properties of $A$ and $\mathcal{A}$ including the null-space property (NSP) and the restricted isometry property (RIP). It ensures the success of the low-rank matrix completion algorithms presented in Sect. 2, restricted isometry constants are introduced in Definition 3.2 and Definition 3.3 , the success of sparse vectors recovery and of low-rank matrices completion are then established under some conditions on these constants for our models in (1.6), (1.7).
\subsection{Recovery Guarantees for Compressed Sensing}
\begin{defn}
A matrix $A\in\mathbb{R}^{n\times p}$ is said to satisfies the null-space property relative to a set $\mathcal{S}\subset[p]$ if
$$\|h_{\mathcal{S}}\|_{1}\leq\|h_{\bar{\mathcal{S}}}\|_{1},\eqno (3.1)$$
for all $h\in NULL(A)/\{0\}$ \cite{R. Gribonval{2003},D. Donoho{2001}}.
\end{defn}
It is said to satisfy the null-space property of order $k$ if it satisfies the null-space property relative to any set $\mathcal{S}\subset[p]$ with $card(\mathcal{S})\leq k$. Given every vector $x\in\mathbb{R}^{p}$ supported on a set $\mathcal{S}$ is the unique solution of (1.3) if and only if $A$ satisfies the null-space property relative to $\mathcal{S}$. Then, we extend the necessary and sufficient NSP condition to our augment model (1.6).

\begin{thm}(NSP condition).

We choose the augmented regularization term $\Phi$ introduced in (1.5). Problem (1.6) uniquely recovers k-sparse vector $x_{0}$ from measurement $Ax_{0}=y$ if
$$(1+\beta\cdot\alpha)\|h_{\mathcal{S}}\|_{1}\leq\|h_{\bar{\mathcal{S}}}\|_{1}\eqno (3.2)$$
hold for all vectors $h\in NULL(A)$ and coordinate sets $\mathcal{S}$ of cardinality $|\mathcal{S}|\leq k$.
\end{thm}

\begin{pf}
$\|x_{0}+h\|_{1}+\beta\Phi(x_{0}+h)$\\
$=\|x_{0}+h_{\mathcal{S}}\|_{1}+\beta\Phi(x_{0}+h_{\mathcal{S}})+\|h_{\bar{\mathcal{S}}}\|_{1}+\beta\Phi(h_{\bar{\mathcal{S}}})$\\
$\geq\|x_{0}\|_{1}-\|h_{\mathcal{S}}\|_{1}+\beta\Phi(x_{0})+\beta\Phi(x_{0}+h_{\mathcal{S}})-\beta\Phi(x_{0})+\|h_{\bar{\mathcal{S}}}\|_{1}+\beta\Phi(h_{\bar{\mathcal{S}}})$\\
$\geq[\|x_{0}\|_{1}+\beta\Phi(x_{0})]-\|h_{\mathcal{S}}\|_{1}+\beta\Sigma_{i\in\mathcal{S}}\frac{\alpha(|x_{i}+h_{i}|-|x_{i}|)}{(1+\alpha|x_{i}+h_{i}|)\cdot(1+\alpha|x_{i}|)}+\alpha\|h_{\bar{\mathcal{S}}}\|_{1}+\Psi(h_{\bar{\mathcal{S}}})$\\
$\geq[\|x_{0}\|_{1}+\beta\Phi(x_{0})]+[\|h_{\bar{\mathcal{S}}}\|_{1}-(1+\beta\cdot\alpha)\|h_{\mathcal{S}}\|_{1}]+\beta\Phi(h_{\bar{\mathcal{S}}})$,\\
where the first inequality from the triangle inequality and the second follows from $\frac{\alpha(|x_{i}+h_{i}|-|x_{i}|)}{(1+\alpha|x_{i}+h_{i}|)\cdot(1+\alpha|x_{i}|)}\geq-\alpha|h_{i}|$. Since $\Phi(h_{\bar{\mathcal{S}}})>0$, and $\|x_{0}+h\|_{1}+\beta\Phi(x_{0}+h)$ is strictly larger than $\|x_{0}\|_{1}+\beta\Phi(x_{0})$, so we can derive inequality (3.2).
\end{pf}

\begin{defn}
The kth restricted isometry constant $\delta_{k}=\delta_{k}(A)$ of matrix $A\in\mathbb{R}^{n\times p}$ is the smallest $\delta\geq0$ such that
$$(1-\delta)\|x\|_{2}^{2}\leq\|Ax\|_{2}^{2}\leq(1+\delta)\|x\|_{2}^{2},\eqno (3.3)$$
for all k-sparse vectors $x\in\mathbb{R}^{p}$ \cite{T. Tao{2005}}.
\end{defn}
We say that $A$ satisfies the restricted isometry property if $\delta_{k}$ is small for reasonably large $k$, then we establish the success of sparse recovery via augment model (1.6) for measurement matrices with small restricted isometry constants.
\begin{thm}
Assume that $x_{0}\in \mathbb{R}^{n}$ is k-sparse. If $A$ satisfies RIP with $\delta_{2k}\leq0.4663$ and $\beta\leq\frac{1}{20\alpha}$, then $x_{0}$ is the unique minimizer of (1.6) given by measurement $Ax_{0}=y$.

\end{thm}
\begin{pf}
\cite{Q. Mo{2011}} shows that any vectors $h\in NULL(A)$ satisfies
$$\|h_{\mathcal{S}}\|_{1}\leq\theta_{2k}\|h_{\bar{\mathcal{S}}}\|_{1},\eqno (3.4)$$
where
$$\theta_{2k}:=\sqrt{\frac{4(1+5\delta_{2k}-4\delta^{2}_{2k})}{(1-\delta_{2k})(32-25\delta_{2k})}},\eqno (3.5)$$
From (3.3), we have
$$\beta\leq\frac{1}{\alpha}(\frac{1}{\theta_{2k}}-1),\eqno (3.6)$$
for $\delta_{2k}=0.4663$, we obtain $\frac{1}{\theta_{2k}}-1\approx0.05\geq\alpha\cdot\beta$.
\end{pf}
Remark1: For (1.3) to recover any k-sparse vector uniformly, \cite{E.J. Candes{2008}} shows the sufficiency of $\delta_{2k}<0.4142$ and improved to $\delta_{2k}<0.4404$ \cite{Ming-Jun Lai{2013}}, $\delta_{2k}<0.4531$ \cite{S. Foucart{2009}}, $\delta_{2k}<0.4652$ \cite{S. Foucart{2010}}, $\delta_{2k}<0.4931$ \cite{Q. Mo{2011}} and the bound is still being improved. \\
Remark2: In general, we choose $\alpha=0.5$ in PF(1.5), so we have $\beta\leq0.1$.

Next, it shows that the condition $\delta_{2k}\leq0.4378$ is actually sufficient to guarantee stable recovery of $x$ via augmented model (2.1).
\begin{thm}
Let $x_{0}\in\mathbb{R}^{n}$ be a arbitrary vector, $\mathcal{S}$ be the coordinate set of its k largest components in magnitude. Let $x^{*}$ be the solution of and error vector $h=x^{*}-x_{0}$ satisfy
$$\|h_{\bar{\mathcal{S}}}\|\leq C_{1}\|h_{\mathcal{S}}\|+C_{2}\|(x_{0})_{\bar{\mathcal{S}}}\|_{1},\eqno (3.7)$$
where
$$C_{1}=\frac{1+\alpha\cdot\beta}{1-\alpha\cdot\beta}~,~C_{2}=\frac{2}{1-\alpha\cdot\beta}.\eqno (3.8)$$

\end{thm}

\begin{pf}
Since $x^{*}=x_{0}+h$ is the minimizer of (1.6), we have
$$\|x_{0}+h\|_{1}+\beta\Phi(x_{0}+h)\leq\|x_{0}\|_{1}+\beta\Phi(x_{0}).\eqno (3.9)$$
We have
$\|x_{0}+h\|_{1}+\beta\Phi(x_{0}+h)$\\
$=\|(x_{0})_{\mathcal{S}}+h_{\mathcal{S}}\|_{1}+\beta\Phi((x_{0})_{\mathcal{S}}+h_{\mathcal{S}})+\|(x_{0})_{\bar{\mathcal{S}}}+h_{\bar{\mathcal{S}}}\|_{1}+\beta\Phi((x_{0})_{\bar{\mathcal{S}}}+h_{\bar{\mathcal{S}}})$\\
$\geq\|(x_{0})_{\mathcal{S}}\|_{1}-\|h_{\mathcal{S}}\|_{1}+\beta\Phi((x_{0})_{\mathcal{S}})-\|(x_{0})_{\bar{\mathcal{S}}}\|_{1}+\|h_{\bar{\mathcal{S}}}\|_{1}+\beta\Phi((x_{0})_{\bar{\mathcal{S}}})$ \\ $+\beta(\Phi((x_{0})_{\mathcal{S}}+h_{\mathcal{S}})-\Phi((x_{0})_{\mathcal{S}}))+\beta(\Phi((x_{0})_{\bar{\mathcal{S}}}+h_{\bar{\mathcal{S}}})-\Phi((x_{0})_{\bar{\mathcal{S}}}))$\\
$\geq[\|x_{0}\|_{1}+\beta\Phi(x_{0})]-2\|(x_{0})_{\bar{\mathcal{S}}}\|_{1}-(1+\alpha\cdot\beta)\|h_{\mathcal{S}}\|_{1}+(1-\alpha\cdot\beta)\|h_{\bar{\mathcal{S}}}\|_{1}$.\\
From (3.9), we have
$$\|h_{\bar{\mathcal{S}}}\|_{1}\leq\frac{1+\alpha\cdot\beta}{1-\alpha\cdot\beta}\|h_{\mathcal{S}}\|+\frac{2}{1-\alpha\cdot\beta}\|(x_{0})_{\bar{\mathcal{S}}}\|_{1}.\eqno (3.10)$$

\end{pf}
\begin{thm}
(see \cite{Ming-Jun Lai{2013}})Let $y=Ax+n$, where $n$ is a arbitrary noise vector with $\|n\|_{2}\leq\epsilon$. If $A$ satisfied RIP with $\delta_{2k}\leq0.4378$, then the solution $x^{*}$ of (2.1) satisfies
$$\|x^{*}-x_{0}\|_{1}\leq C_{3}\cdot\sqrt{k}\|n\|_{2}+C_{4}\cdot\|(x_{0})_{\bar{\mathcal{S}}}\|_{1},\eqno (3.11)$$
$$\|x^{*}-x_{0}\|_{2}\leq C_{5}\cdot\|n\|_{2}+C_{6}\cdot\|(x_{0})_{\bar{\mathcal{S}}}\|_{1}/\sqrt{k},\eqno (3.12)$$
where
$$C_{3}=\frac{2\sqrt{2}(1+C_{1})}{\sqrt{1-\delta_{2k}}(1-C_{1}\theta_{2k})}~,~C_{4}=\frac{(1+\theta_{2k})C_{2}}{1-C_{1}\theta_{2k}},\eqno (3.13)$$
and
$$C_{5}=\frac{2}{\sqrt{1-\delta_{2k}}}\{\frac{4C_{1}}{1-C_{1}\theta_{2k}}\sqrt{\frac{2-\delta_{2k}}{(1-\delta_{2k})(32-25\delta_{2k})}}+1\},\eqno (3.14)$$
$$C_{6}=\frac{2C_{2}}{1-C_{1}\theta_{2k}}\sqrt{\frac{2(2-\delta_{2k})}{(1-\delta_{2k})(32-25\delta_{2k})}}.\eqno (3.15)$$
\end{thm}
\subsection{Recovery Guarantees for Matrix Recovery }
It's easy to extended the NSP and RIP condition to low-rank matrix recovery, first, let us introduce some definitions and properties. $\|X\|_{*}=\sum_{i=1}^{r}\sigma_{i}(X)$, $\|X\|_{F}=\sqrt{\sum_{i=1}^{r}\sigma^{2}_{i}(X)}$ denote the unclear and Frobenius norm of $X$ respectively, where $\sigma_{i}$ is the $i$ largest singular value of $X$ and $r$ is the number of singular value.

Let $X$ and $W$ be two matrices of the same size, we have $\sum^{p}_{i=r}\varphi(\sigma_{i}(X)-\sigma_{i}(W))\leq\hat{\Phi}(X-W)$, because  $|\sigma_{i}(X)-\sigma_{i}(W)|\leq|\sigma_{i}(X-W)|$, for $i=1,\ldots,r$ and $\varphi$ is a increasing function.
\begin{thm}
Problem (1.7) uniquely recovers all matrices $X$ of rank $r$ or less from measurement $\mathcal{A}X=b$ if
$$(1+\beta\cdot\alpha)\sum_{i=1}^{r}\sigma_{i}(H)\leq\sum_{i=r+1}^{m}\sigma_{i}(H),\eqno (3.16)$$
holds for all matrices $H\in NULL(\mathcal{A})$.
\end{thm}
\begin{pf}
$\|X+H\|_{*}+\beta\hat{\Phi}(X+H)\geq\beta\sum_{i=1}^{m}(\sigma_{i}(X)+\sigma_{i}(H))++\beta\sum_{i=1}^{m}\varphi(\sigma_{i}(X)+\sigma_{i}(H))$\\
$\geq\sum_{i=1}^{r}\sigma_{i}(X)-\sum_{i=1}^{r}\sigma_{i}(H)+\sum_{i=r+1}^{m}\sigma_{i}(H)+\beta\sum_{i=1}^{r}\varphi(\sigma_{i}(X)+\sigma_{i}(H))+\beta\sum_{i=1+r}^{m}\varphi(\sigma_{i}(H))$\\
$\geq[\|X\|_{*}+\beta\hat{\Phi}(X)]+[\sum_{i=r+1}^{m}\sigma_{i}(H)-(1+\beta\cdot\alpha)\sum_{i=1}^{r}\sigma_{i}(H)]+\beta\sum_{i=1+r}^{m}\varphi(\sigma_{i}(H))$.
\end{pf}

\begin{defn}
for a linear map $\mathcal{A}:\mathbb{R}^{n_{1}\times n_{2}}\rightarrow \mathbb{R}^{n_{3}}$ and for $r\leq m=\min\{n_{1},n_{2}\}$, the rank restricted isometry constant $\delta_{r}=\delta_{r}(\mathcal{A})$ is the defined as the smallest $\delta\geq0$ such that
$$(1-\delta)\|X\|_{F}^{2}\leq\|\mathcal{A}(X)\|_{2}^{2}\leq(1+\delta)\|X\|_{F}^{2},\eqno (3.17)$$
for all matrices $X\in\mathbb{R}^{n_{1}\times n_{2}}$ of rank at most $r$ \cite{B. Recht{2010}}.
\end{defn}
\begin{thm}
(RIP condition for exact recovery). Let $X$ be a matrix with rank $r$ or less, the augment model (1.7) exactly recovers $X$ from measurement $b=\mathcal{A}(X)$ if $\mathcal{A}$ satisfies the RIP condition with $\delta_{2r}\leq0.4663$.
\end{thm}
\begin{pf}
In \cite{Ming-Jun Lai{2013}}, establishes that any $H\in NULL(\mathcal{A})$ satisfy $\sum_{i=1}^{r}\sigma_{i}(H)\leq\theta_{2r}\sum_{i=r+1}^{m}\sigma_{i}(H)$, hence (3.16) holds if $(1+\alpha\cdot\beta)^{-1}\geq\theta_{2r}$.
\end{pf}
\begin{thm}
(RIP condition for stable recovery) Let $X\in\mathbb{R}^{n_{1}\times n_{2}}$ be an arbitrary matric, and let $b=\mathcal{A}X+n$, where $\mathcal{A}$ is a linear operator and $n$ is an arbitrary noise. If $\mathcal{A}$ satisfies the RIP with $\delta_{2r}\leq0.4378$, then, the solution $X^{*}$ of (2.2) satisfies the error bounds
$$\|X^{*}-X\|_{*}\leq\bar{C}_{3}\cdot\sqrt{r}\|n\|_{2}+\bar{C}_{4}\cdot\sum_{i=r+1}^{m}\sigma_{i}(X),\eqno (3.18)$$
$$\|X^{*}-X\|_{F}\leq\bar{C}_{5}\cdot\|n\|_{2}+\bar{C}_{6}\cdot\sum_{i=r+1}^{m}\sigma_{i}(X)/\sqrt{r},\eqno (3.19)$$
$\bar{C}_{3}$, $\bar{C}_{4}$, $\bar{C}_{5}$ and $\bar{C}_{6}$ are given formulas (3.13)-(3.15) in which $\theta_{2k}$ shall be replaced by $\theta_{2r}$.
\end{thm}

\section{Numerical Experiments}
\subsection{Test on Compressed Sensing}
In this subsection, we perform experiments on synthetic data to illustrate the behavior of the augmented nonconvex method and Lasso. The
support $S$ of $x$ is equal to $\{1,\ldots,k\}$, where $k$ is the size of the support. For $i$ in the support of
$x$, $x_i$ is independently drawn from a Gaussian distribution with zero mean and standard deviation $\sigma=1$. The $A_{i}$ are
drawn from a multivariate Gaussian with mean zero and covariance matrix $\Sigma$, where $A_{i}$ is the $i$ column of ensemble $A$. For the first setting, $\Sigma$
is set to the identity, for the second setting, $\Sigma$ is block diagonal with blocks equal to $0.2I+0.8\mathbf{1}\mathbf{1}^{*}$ \cite{E. Grave{2011}}. We perform the experiments $(p=512, n=128)$ for which we report the estimation relative error, which defines as
$$RelErr=\frac{\|x^{opt}-x\|_{2}}{\|x\|_{2}}.$$

The recovery is performed via the augment nonconvex method algorithm, and we use
$$\|x^{(s)}-x^{(s-1)}\|_{2}/\|x^{(s-1)}\|_{2}<10^{-4},$$
and the maximum iteration step $maxit=500$ as stopping criterion.
In Fig 1. we observe that the Lasso performs as well as the augmented nonconvex method with parameter $\alpha=0.5$, $\beta=0.1$ and $\alpha=0.1$, $\beta=0.5$ on very sparse case. But, when the support of $x$ is large, the augmented nonconvex method perform well than Lasso on both two setting \cite{R. Tibshirani{1996}}.
\begin{figure*}[htbp]
\begin{center}
\begin{tabular}{c}
\subfigure[]{\includegraphics[width=0.3\textwidth,height=0.3\textwidth]{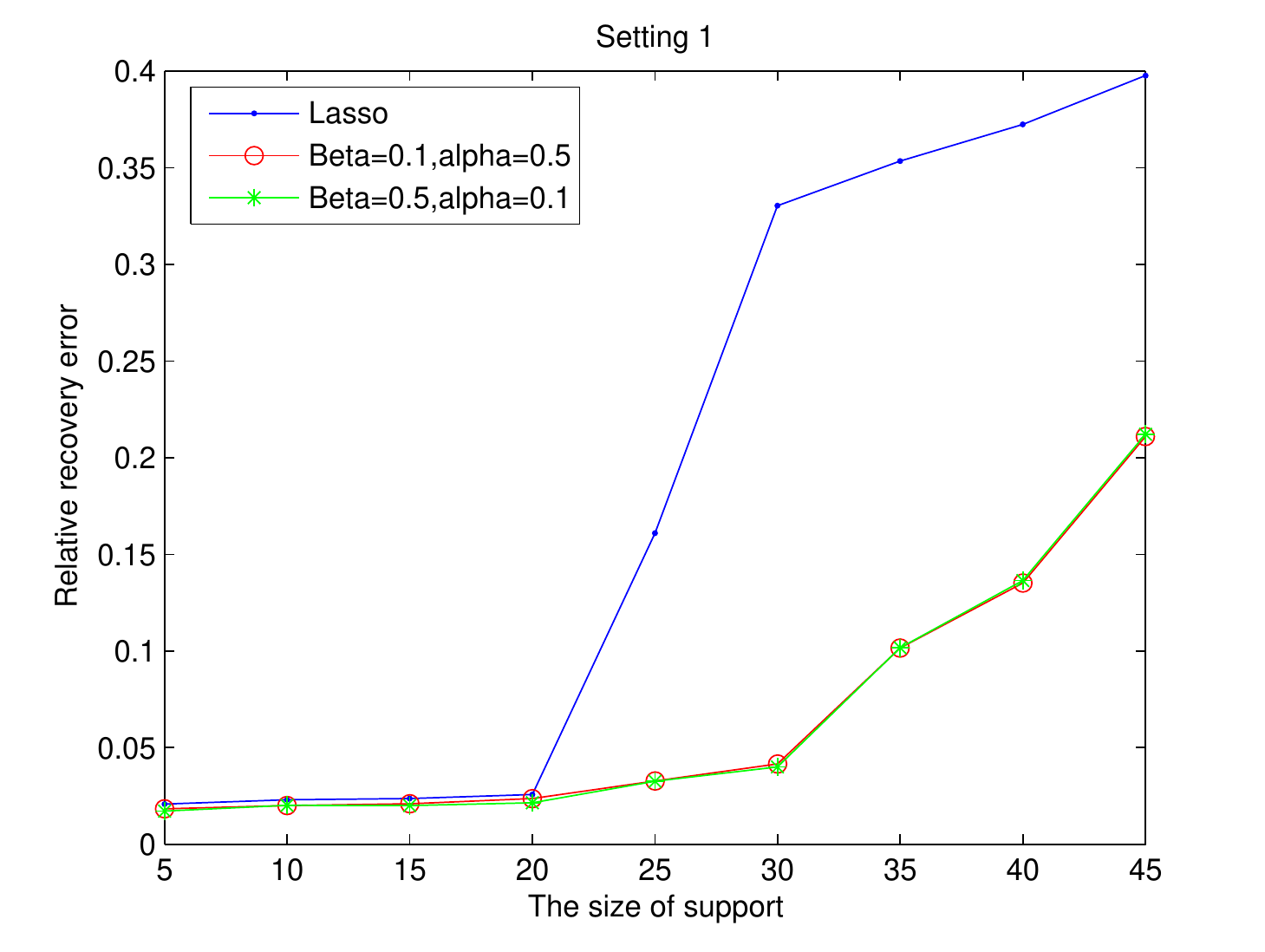}}
\subfigure[]{\includegraphics[width=0.3\textwidth,height=0.3\textwidth]{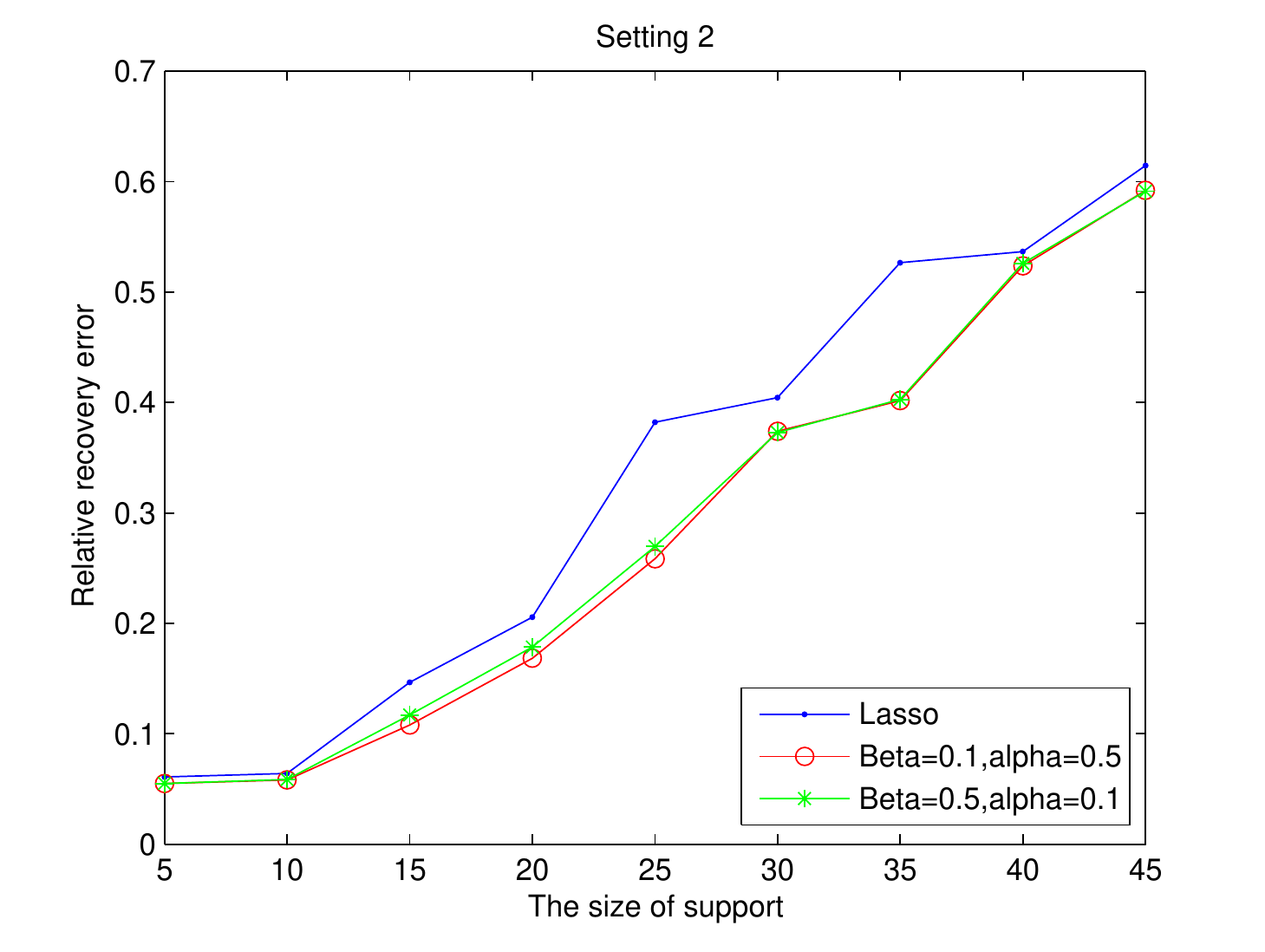}}
\end{tabular}
\end{center}
\caption{(a) Setting 1, the column of ensemble $A_{i}$ with covariance matrix $\Sigma=I$; (b) Setting 2, the column of ensemble $A_{i}$ with covariance matrix $\Sigma=0.2I+0.8\mathbf{1}\mathbf{1}^{*}$.}
\label{f5b}
\end{figure*}
\subsection{Reconstruction from Sparse Fourier Measurement}
In this subsection, we consider the problem of image reconstruction from a limited number of Fourier measurements. In this setting, the operator of (1.1) corresponds to $A=MF$, where $F$ denotes the Fourier transform and $M$ is a masking operator the retains only a subset of the available Fourier coefficients \cite{Stamatios Lefkimmiatis{2015}}, and we use the augmented nonconvex method solve the following problem
$$\min_{f\in\mathbb{R}^{n_{1}\times n_2}}
~\frac{1}{2}\|Af-g\|_{2}^{2}+\lambda\cdot(\|f\|_{TV}+\beta\sum_{t_1,t_2}\varphi(\sqrt{|D_{1}f(t_1,t_2)|^{2}+|D_{2}f(t_1,t_2)|^{2}})),$$
where $\|f\|_{TV}$ is the total variation norm, for discrete $f(t_1,t_2)$, $0\leq t_1\leq n_1,0\leq t_2\leq n_2$ and $D_1$ is the finite difference $D_1f(t_1,t_2)=f(t_1,t_2)-f(t_1-1,t_2)$ and $D_2$ is the finite difference $D_1f(t_1,t_2)=f(t_1,t_2)-f(t_1,t_2-1)$.
\begin{figure*}[htbp]
\begin{center}
\begin{tabular}{c}
\subfigure[]{\includegraphics[width=0.1925\textwidth,height=0.1925\textwidth]{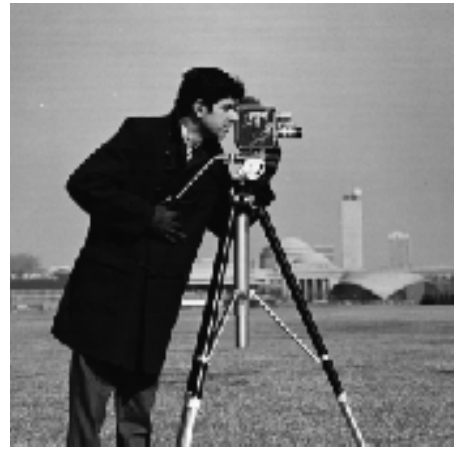}}
\subfigure[]{\includegraphics[width=0.1925\textwidth,height=0.1925\textwidth]{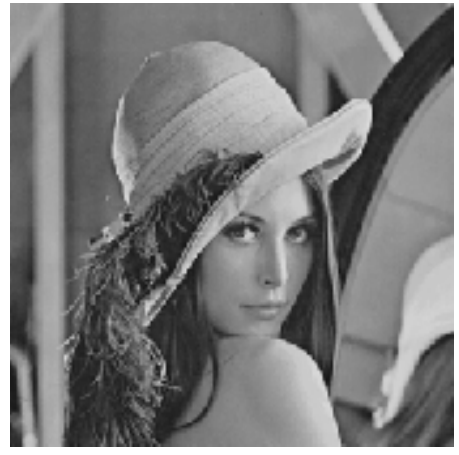}}
\subfigure[]{\includegraphics[width=0.1925\textwidth,height=0.1925\textwidth]{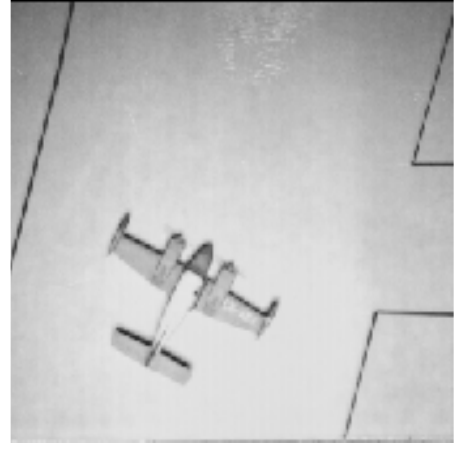}}
\subfigure[]{\includegraphics[width=0.1925\textwidth,height=0.1925\textwidth]{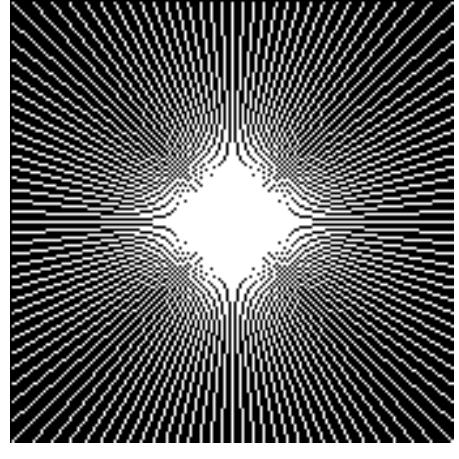}}
\subfigure[]{\includegraphics[width=0.1925\textwidth,height=0.1925\textwidth]{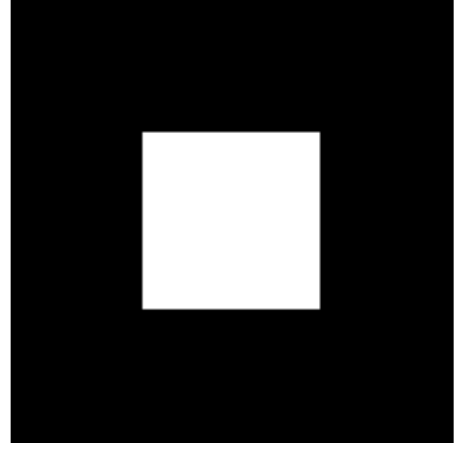}}
\end{tabular}
\end{center}
\caption{(a)-(c) Images with size $256\times256$ and downsample factor=1.5; (d) Radial sampling mask with 64 lines; (e) low-frequency sampling with 40\% portion.}
\label{f5b}
\end{figure*}
The reported experiments are conducted on images shows in Fig. 2. To create the measured data we use two different Fourier sampling patterns, namely, a radial mask with 48, 64 and 80 radial lines and a low-frequency sampling with 30\%, 40\% and 50\% portion. As an additional degradation factor we consider the presence of complex Gaussian noise in Fourier domain of four different levels. These correspond to a signal-to-noise-ratio (SNR) of $\{15, 20, 30, \infty\}$dB, and the last SNR value indicates the absence of noise in the measurements. Peak-signal-to-noise ratio(PSNRs) is used to measure the quality of the restored images, which are defined as
$$PSNR=10\cdot\log\frac{255^{2}}{MSE}[dB]$$
where MSE is the Mean-Squared-Error per pixel. In Table. 1 and Fig. 3 shows that the augmented nonconvex method perform well than TV.
\begin{center}
\begin{table}
\caption{PSNR comparisons on Fourier image reconstruction for several sampling patterned and noise conditions (N-TV represent our proposed nonconvex regularization method).}
\resizebox{\textwidth}{!}{

\begin{tabular}{cr|rrrr|rrrr|rrrr|rrrr|rrrr|rrrr}
\hline
\hline
\multicolumn{ 2}{c}{ Sampling} & \multicolumn{ 4}{|c}{Radial-48 lines} & \multicolumn{ 4}{|c}{Radial-64 lines} & \multicolumn{ 4}{|c}{Radial-80 lines} & \multicolumn{ 4}{|c}{Low-frequency portion  30£¥} & \multicolumn{ 4}{|c}{Low-frequency portion  40£¥} & \multicolumn{ 4}{|c}{Low-frequency portion  50£¥} \\

\multicolumn{ 2}{c|}{ PSNR   } & 15dB & 20dB & 30dB & $\infty$dB & 15dB & 20dB & 30dB & $\infty$dB  & 15dB & 20dB & 30dB & $\infty$dB  & 15dB & 20dB & 30dB & $\infty$dB  & 15dB & 20dB & 30dB & $\infty$dB  & 15dB & 20dB & 30dB & $\infty$dB  \\
\hline
\multicolumn{ 1}{c|}{Lenna} & N-TV & \textbf{38.324} & \textbf{38.659} & \textbf{38.776} & \textbf{38.797} & \textbf{38.961} & \textbf{39.496} & \textbf{39.714} & \textbf{39.741} & \textbf{39.373} & \textbf{40.121} & \textbf{40.441} & \textbf{40.487} & \textbf{38.425} & \textbf{38.549} & \textbf{38.615} & \textbf{38.610} & \textbf{39.024} & \textbf{39.298} & \textbf{39.381} & \textbf{39.406} & \textbf{39.502} & \textbf{39.930} & \textbf{40.112} & \textbf{40.141} \\
\hline
\multicolumn{ 1}{c|}{} &   TV & 38.166 & 38.367 & 38.475 & 38.497 & 38.885 & 39.156 & 39.344 & 39.371 & 39.331 & 39.776 & 40.054 & 40.092 & 38.230 & 38.362 & 38.429 & 38.448 & 38.855 & 39.106 & 39.222 & 39.246 & 39.340 & 39.732 & 39.944 & 39.978 \\
\hline
\multicolumn{ 1}{c|}{Cameraman} & N-TV & \textbf{37.799} & \textbf{37.905} & \textbf{37.973} & \textbf{37.983} & \textbf{38.276} & \textbf{38.512} & \textbf{38.598} & \textbf{38.591} & \textbf{38.758} & \textbf{38.997} & \textbf{39.107} & \textbf{39.122} & \textbf{37.070} & \textbf{37.129} & \textbf{37.142} & \textbf{37.144} & \textbf{37.629} & \textbf{37.777} & \textbf{37.838} & \textbf{37.839} & \textbf{38.108} & \textbf{38.383} & \textbf{38.469} & \textbf{38.481} \\
\hline
\multicolumn{ 1}{c|}{} &   TV & 37.441 & 37.571 & 37.668 & 37.637 & 37.927 & 38.206 & 38.328 & 38.350 & 38.544 & 38.798 & 38.993 & 39.018 & 36.899 & 36.973 & 37.025 & 37.005 & 37.523 & 37.642 & 37.712 & 37.732 & 38.037 & 38.238 & 38.362 & 38.369 \\
\hline
\multicolumn{ 1}{c|}{Airplane} & N-TV & \textbf{39.555} & \textbf{40.426} & \textbf{40.782} & \textbf{40.819} & \textbf{40.065} & \textbf{40.780} & 41.103 & 41.148 & 39.751 & \textbf{41.496} & \textbf{42.380} & \textbf{42.482} & \textbf{38.426} & \textbf{38.583} & \textbf{38.645} & \textbf{38.639} & \textbf{38.847} & \textbf{39.174} & \textbf{39.257} & \textbf{39.268} & \textbf{39.078} & \textbf{39.687} & \textbf{39.839} & \textbf{39.853} \\
\hline
\multicolumn{ 1}{c|}{} &   TV & 39.236 & 39.893 & 40.401 & 40.509 & 39.718 & 40.571 & \textbf{41.486} & \textbf{41.679} & \textbf{40.039} & 40.934 & 42.164 & 42.418 & 38.257 & 38.384 & 38.490 & 38.456 & 38.696 & 38.975 & 39.116 & 39.112 & 38.960 & 39.476 & 39.472 & 39.765 \\
\hline
\hline
\end{tabular}
}
\end{table}
\end{center}
\begin{figure*}[htbp]
\begin{center}
\begin{tabular}{c}
\subfigure[]{\includegraphics[width=0.1925\textwidth,height=0.1925\textwidth]{camera.pdf}}
\subfigure[]{\includegraphics[width=0.1925\textwidth,height=0.1925\textwidth]{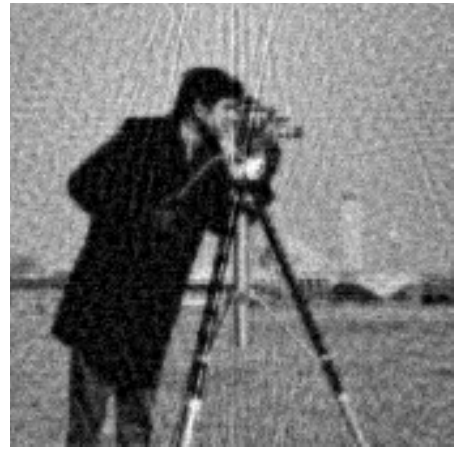}}
\subfigure[]{\includegraphics[width=0.1925\textwidth,height=0.1925\textwidth]{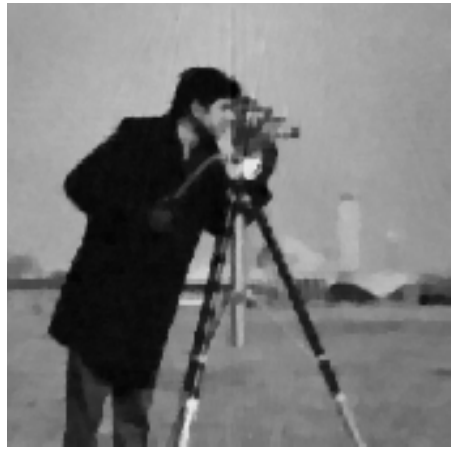}}
\subfigure[]{\includegraphics[width=0.1925\textwidth,height=0.1925\textwidth]{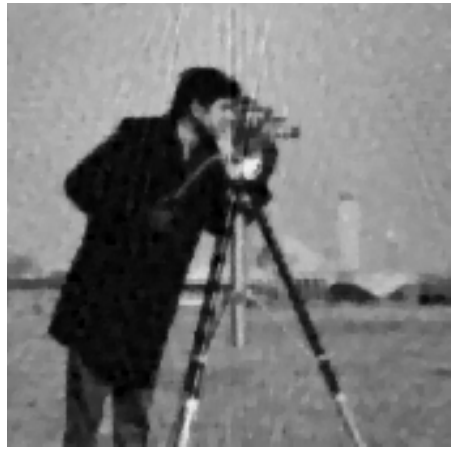}}
\end{tabular}
\end{center}
\caption{Reconstruct the image from Fourier data sampled with 64 radial lines and 20-dB SNR. (a) Original image; (b) Back-projiect image; (c) Our proposed  nonconvex regularization method; (d) TV reconstruction.}
\label{f5b}
\end{figure*}
\subsection{Test on Matrix Completion}
In our numerical experiments, $n_{1}$ and $n_2$ represent the matrix dimension, $r$ is the rank of original matrix, and $n_3$ denotes the number of measurement. Given $r\leq min(n_1,n_2)$, we generate $X=X_{L}X^{*}_{R}$, where matric $X_{L}\in\mathbb{R}^{n_{1}\times r}$ and $X_{R}\in\mathbb{R}^{n_{2}\times r}$ are generated with independent identically distributed Gaussian entries. The subset $\Omega$ of $n_3$ elements is selected uniformly at random entries from $\{(i,j):i=1,\ldots,n_1,j=1,\ldots,n_2\}$ \cite{Z.F. Jin{2015}}. The linear measurement $b$ are set to be $b=\mathcal{A}(X)+n$, where $n$ is the additive Gaussian noise of zero mean and standard deviation $\sigma$, which will be specified in different test data sets. We use $sr=n_{3}/(n_1n_2)$ to denote the sampling ratio, and $dr=r(n_1+n_2-r)$ to denote the number of degree of freedom for a real-valued rank $r$ matrix. As mentioned in \cite{M. Malek-Mohammadi{2014}}, when the ratio $n_{3}/dr$ is greater than 3, the problem can be viewed as an easy problem. On the contrary, it is called as a hard problem.

In this subsection, we apply the proposed augmented nonconvex method for solving the matrix completion problem (2.4). In order to illustrate the performance of this method, we compare the augmented nonconvex method with the nuclear-norm model \cite{B. Recht{2009}} and the augmented Nuclear-Norm model with $\alpha=50$ \cite{Ming-Jun Lai{2013}}.

The recovery is performed via the augment nonconvex method algorithm, and we use
$$\|X^{(s)}-X^{(s-1)}\|_{F}/\|X^{(s-1)}\|_{F}<10^{-8},$$
and the maximum iteration step $maxit=2000$ as stopping criterion. Our computational results are displayed in Table 2. We choose $n_1=n_2$, noise level $\sigma=1e-3$, and the relative error of the reconstruction matrix $X$ is
$$RelErr=\frac{\|X^{opt}-X\|_{F}}{\|X\|_{F}},$$
and it shows that the augment nonconvex method (the last column) can get higher accuracy than others.
\begin{center}
\begin{table}[h]
\caption{Numerical result of nuclear-norm model (Nuclear) \cite{B. Recht{2009}}, augmented nuclear-norm model (Aug-Nuclear) \cite{Ming-Jun Lai{2013}} and our proposed nonconvex regularizaion method (N-Nuclear) for matrix completion problems.}
\centering
\begin{tabular}{r|r|rrr}
\hline
\hline
\multicolumn{ 1}{c|}{$(n_1,r)$} & \multicolumn{ 1}{|c|}{$n_{3}/dr$} & Nuclear & Aug-Nuclear & N-Nuclear \\

\multicolumn{ 1}{c|}{} & \multicolumn{ 1}{|c|}{} & RelErr & RelErr & RelErr \\
\hline
(100,10) & 2.632& 8.01e-04 & 9.30e-04 & 9.48e-05 \\

(200,20) & 2.632 & 9.02e-04 & 9.71e-04 & 5.78e-05 \\

(300,30) & 2.632 & 7.88e-04 & 4.35e-04 & 4.50e-05 \\

(400,40) & 2.632 & 6.63e-04 & 4.90e-04 & 5.29e-05 \\

(500,50) & 2.632& 6.57e-04 & 5.25e-04 & 5.23e-05 \\
\hline
\hline
\end{tabular}
\end{table}
\end{center}

Finally, we test the augmented nonconvex method for recovering two real corrupted gray scale image. at first, we use SVD to obtain the low-rank-50 images. Then we randomly select $40\%$ samples from the low-rank image, which corrupted image with noise level $\sigma=1e-3$. Finally, these corrupted images are corrupted images are recovered by the proposed nonconvex regularization method and the nuclear-norm model. From Fig. 1, it showed that the quality of image (c) restored by augmented nonconvex method is better than the image (d) restored by nuclear-norm model.

\begin{figure*}[htbp]
\begin{center}
\begin{tabular}{c}
\subfigure{\includegraphics[width=0.1925\textwidth,height=0.1925\textwidth]{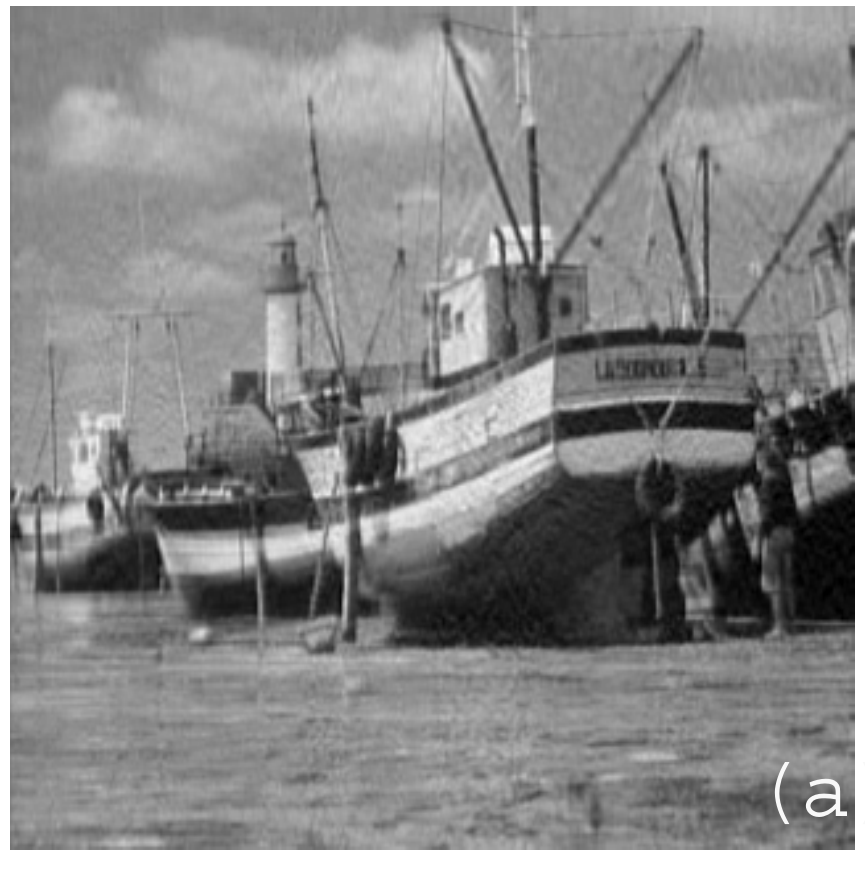}}
\subfigure{\includegraphics[width=0.1925\textwidth,height=0.1925\textwidth]{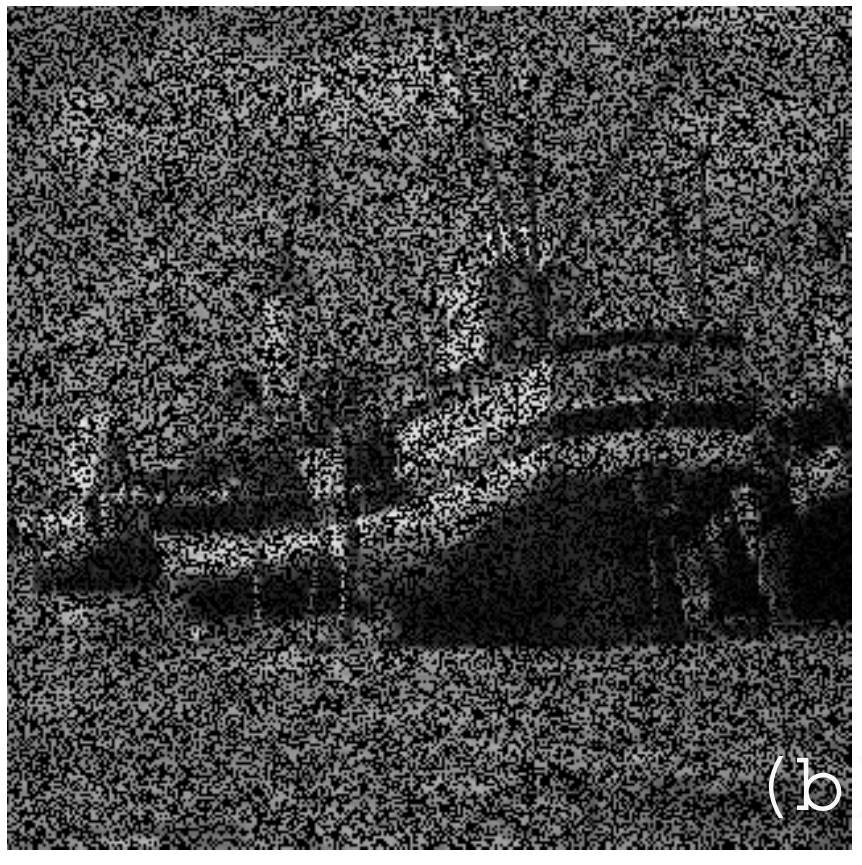}}
\subfigure{\includegraphics[width=0.1925\textwidth,height=0.1925\textwidth]{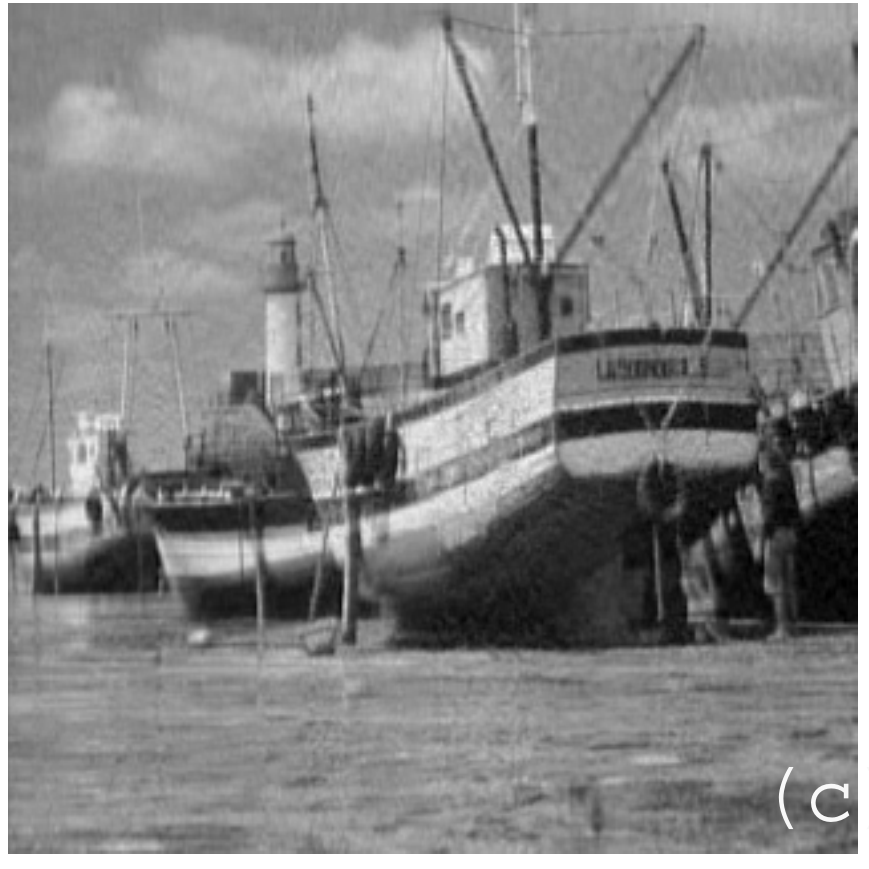}}
\subfigure{\includegraphics[width=0.1925\textwidth,height=0.1925\textwidth]{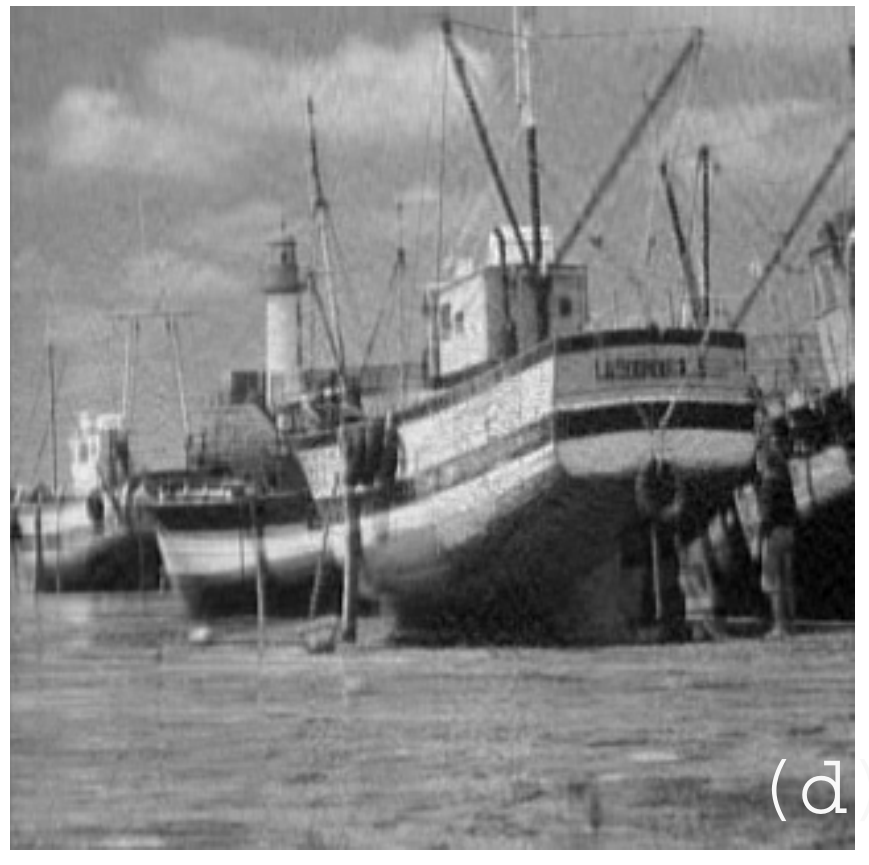}}\\
\subfigure{\includegraphics[width=0.1925\textwidth,height=0.1925\textwidth]{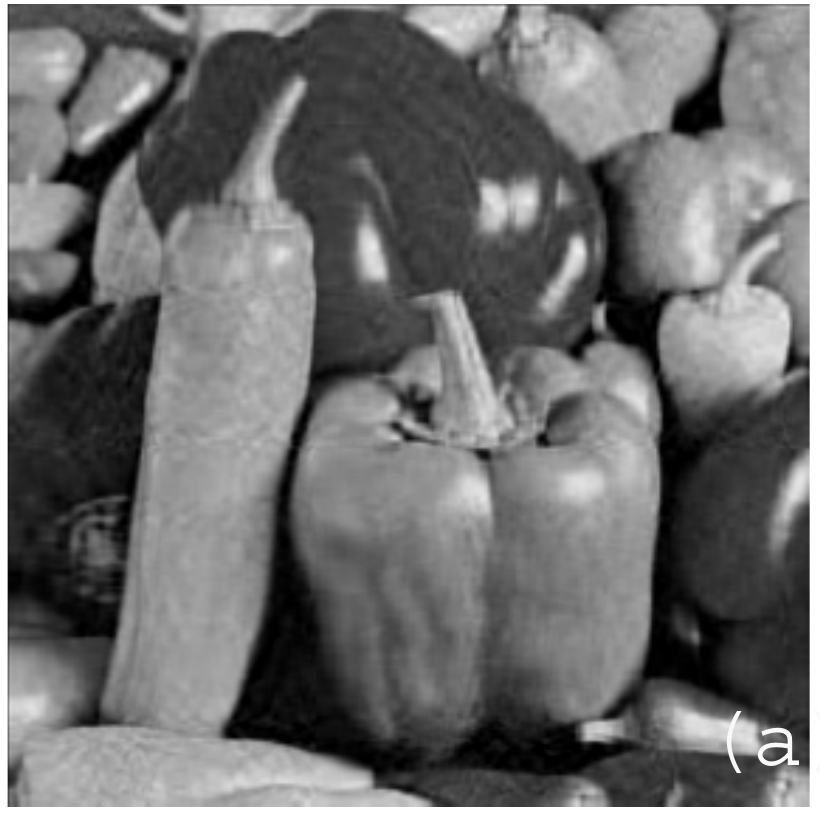}}
\subfigure{\includegraphics[width=0.1925\textwidth,height=0.1925\textwidth]{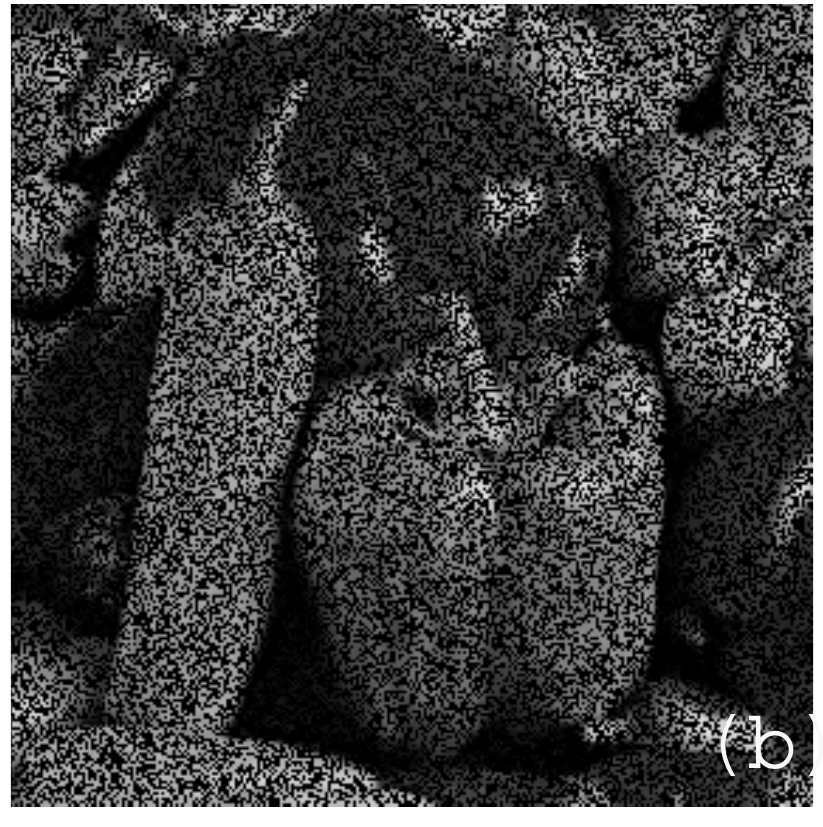}}
\subfigure{\includegraphics[width=0.1925\textwidth,height=0.1925\textwidth]{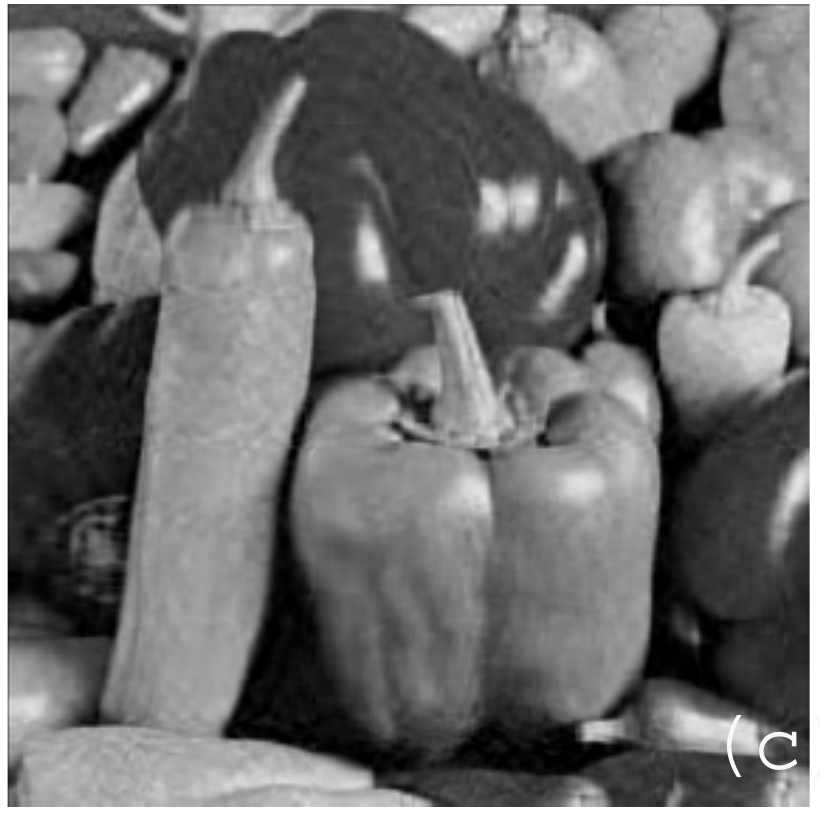}}
\subfigure{\includegraphics[width=0.1925\textwidth,height=0.1925\textwidth]{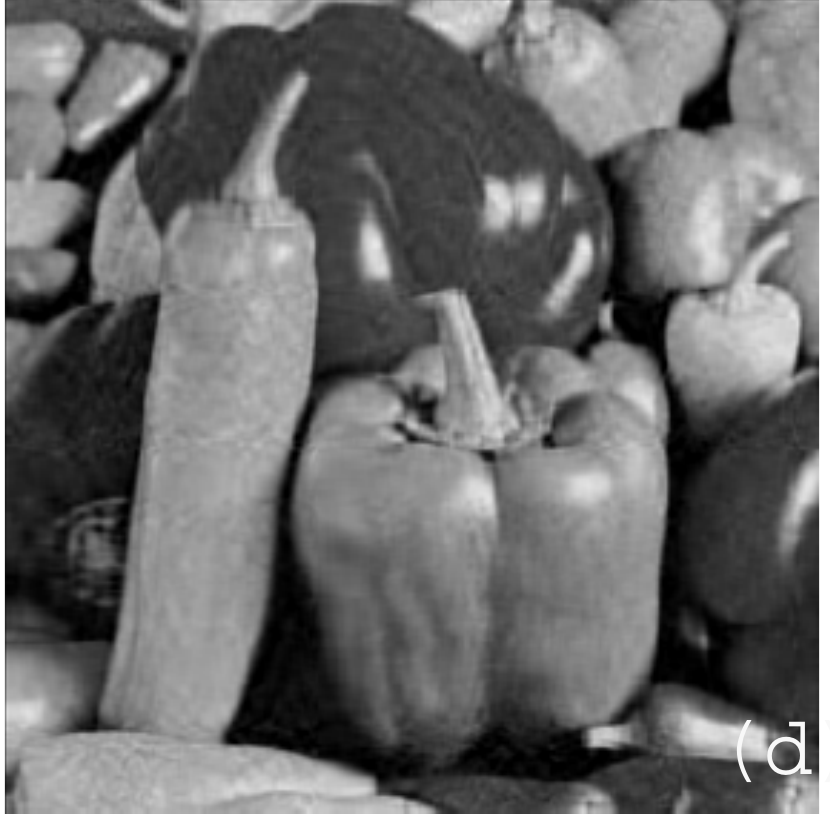}}
\end{tabular}
\end{center}
\caption{(a) Corresponding low rank image with $n_1=n_2=512$, $r=50$ and downsamplefactor=1.8; (b) Randomly masked images from rank 50 with $sr=40\%$, $\sigma=1e-3$; (c) Recovered images by our proposed nonconvex regularization method [ErrRel=7.93e-03 (first image), 1.01e-02 (second image)]; (d) Recovered images by nuclear-norm model [ErrRel=2.04e-02 (first image), 2.36e-02 (second image)].}
\label{f5b}
\end{figure*}

\section{conclusions}
In this paper, we given the augmented model, and introduced the nonconvex and nonsmooth penalty function for low-rank matrix completion and sparse vector recovery. Then, we developed the alternating minimization scheme for solving the proposed problem and give the convergence result of the proposed method. In addition, we showed the recovery guarantees for augmented model of compressed sensing and low-rank matrix completion respectively, including NSP and RIP.
At last, some numerical experiment results of our augmented model have been showed on simulated and real data and performs well. However, the unclear norm measures the low-rank property of $X$ without considering the interelement of singular value correlations. When the singular values have high correlations, the nuclear norm is known to have stability problems. In the future research work, We desire to measure the low-rank property of $X$ at group level and have
all singular value within a group become nonzero (or zero) simultaneously, and also show the recovery guarantees at group level.
\section*{Acknowledgments}

This research was supported by the National Science Foundation of China under Grant 61179039 and the National Key Basic Research Development Program(973 Program) of China under Grant 2011CB707100.

\section*{Appendix: Algorithm for Sparse Vector Recovery}
\begin{tabular}{l}
 \hline
Algorithm To solve the Minimum Value of (2.1) \\
\hline
Step 1: Initialize $x^{(0)}$ and $s=1$;\\
Step 2: Update $x$ and $w$ until the convergence\\
~~~~~w-step:\\
~~~~~~~$w^{(s)}=\arg\min_{w\in\mathbb{R}^{p}}J(w,x^{(s-1)})$,\\
~~~~~~~$w^{(s)}_{i}=x^{(s-1)}_{i}\max\{1-\frac{\tau}{|x^{(s-1)}_{i}|},0\}$,\\
~~~~~~~~~~~~~~~~~~~~for $i=1,\ldots,p$ and $\tau=\frac{\lambda(1+\alpha)}{\rho}$.\\
~~~~~x-step:\\
~~~~~~~$x^{(s)}=x^{(s-1)}+\Delta x$, where\\
~~~~~~~$(A^{*}A+\rho I)\Delta x=-\nabla J(x,w^{(s)})$.\\
(Here the iteration index is the superscript $s$.)\\
\hline
\end{tabular}\\

\newpage
\section*{References}

\bibliography{mybibfile}

\end{document}